\documentclass[12pt,reqno]{amsart}
\usepackage{amssymb,delarray}
\usepackage{graphicx}

\usepackage{epsf,amsfonts,amscd,latexsym}

\def\leq{\leqslant}
\def\geq{\geqslant}

\def\pr{{\sf pr}}

\def\pr{\mathsf{pr}}

\newtheorem{thm}{Theorem}[section]
\newtheorem{lem}[thm]{Lemma}
\newtheorem{prop}[thm]{Proposition}

\newtheorem{conj}[thm]{Conjecture}
\newtheorem{ex}[thm]{Example}
{\catcode`\@=11
\gdef\n@te#1#2{\leavevmode\vadjust{%
 {\setbox\z@\hbox to\z@{\strut#1}%
  \setbox\z@\hbox{\raise\dp\strutbox\box\z@}\ht\z@=\z@\dp\z@=\z@%
  #2\box\z@}}}
\gdef\leftnote#1{\n@te{\hss#1\quad}{}}
\gdef\rightnote#1{\n@te{\quad\kern-\leftskip#1\hss}{\moveright\hsize}}
\gdef\?{\FN@\qumark}
\gdef\qumark{\ifx\next"\DN@"##1"{\leftnote{\rm##1}}\else
 \DN@{\leftnote{\rm??}}\fi{\rm??}\next@}}

\begin{document}
\baselineskip=14.pt plus 2pt 

\title[
] {On the Alexander polynomials of Hurwitz curves}

 \author[Vik.S.~Kulikov]{Vik.S.~Kulikov}
\address{Steklov Mathematical Institute}
 \email{kulikov@mi.ras.ru}
\dedicatory{} \subjclass{}
\thanks{The author was partially supported  by the RFBR 
(
02-01-00786)}
\keywords{}
\begin{abstract}
Properties of the Alexander polynomials of Hurwitz curves are
investigated. A complete description of the set of the Alexander
polynomials of irreducible Hurwitz curves in the terms of their
roots is given.
\end{abstract}

\maketitle
\setcounter{tocdepth}{2}


\def\st{{\sf st}}

\setcounter{section}{-1}
\section{Introduction}

In \cite{Gr-Ku}, investigation of properties of the Alexander
polynomials of Hurwitz curves in $\mathbb C\mathbb P^2$ were
started. Recall briefly definitions of a Hurwitz curve (with
respect to a linear projection $\text{pr}: \mathbb C\mathbb P^2
\to \mathbb C\mathbb P^1$) and its Alexander polynomial. Let
$\mathbb C^2_i$ be two copies of the affine plane $\mathbb C^2$,
$i=1,2$, with coordinates $(u_i,v_i)$, $u_2=1/u_1$ and
$v_2=v_1/u_1$, which cover $\mathbb C\mathbb P^2\setminus
p_\infty$ (where $p_{\infty}$ is the center of the projection)
 such that $\text{pr}$ is given by $(u_i,v_i)\mapsto u_i$
in the charts $\mathbb C^2_i$. A set $\bar H\subset \mathbb
C\mathbb P^2\setminus \{p_\infty\}$, closed in $\mathbb C\mathbb
P^2$, is called a {\it Hurwitz curve of degree} $m$ if, for
$i=1,2$, $\bar H\cap \mathbb C^2_i$ coincides with the set of
zeros of an equation \[ F_i(u_i,v_i):=v_i^m
+\sum_{j=0}^{m-1}c_{j,i}(u_i)v_i^j=0
\]
such that
\begin{itemize}
\item[($i$)] $F_i (u_i, v_i)$ is a  $C^\infty$-smooth complex valued function in
$\mathbb C^2$;
\item[($ii$)] the function $F_i(u_i,v_i)$ has only
a finite number of critical values, that is, there are finitely
many values of $u_i$, say $u_{i,1},\dots, u_{i,n_i}$, such that
the polynomial equation
\begin{equation} \label{def} v_i^m
+\sum_{j=0}^{m-1}c_{j,i}(u_{i,0})v_i^j=0\end{equation} has no
multiple roots for $u_{i,0}\not\in \{ u_{i,1},\dots, u_{i,n_i}
\}$; \item[($iii$)] if $v_{i,j}$ is a multiple root of equation
(\ref{def}) for $u_{i,j}\in \{ u_{i,1},\dots, u_{i,n_i} \}$, then,
in a neighbourhood of the point $(u_{i,j},v_{i,j})$, the set $\bar
H$ coincides with the solution a complex analytic equation.
\end{itemize}

A Hurwitz curve $\bar H$ is called {\it irreducible} if $\bar
H\setminus M$ is connected for any finite set $M\subset \bar H$,
and we say that a Hurwitz curve $\bar H$ {\it consists of $k$
irreducible components} if $$k=\max \# \{ \text{connected
components of}\, \bar H\setminus M\},$$ where the maximum is taken
over all finite sets $M\subset \bar H$.

Let $H$ be an {\it affine Hurwitz curve}, that is, $H=\bar H\cap
(\mathbb C\mathbb P^2\setminus L_{\infty})$, where $L_{\infty}$ is
a line which is a fibre of $\text{pr}$ being in general position
with respect to $\bar H$. Then the fundamental group
$\pi_1=\pi_1(\mathbb C\mathbb P^2\setminus (\bar H\cup
L_{\infty}))$ does not depend on the choice of $L_{\infty}$ and
belongs to the class $\mathcal C$ of so called $C$-groups.

By definition, a {\it $C$-group} is a group together with  a
finite presentation
\begin{equation} \label{zero}
G_W=\langle x_1,\dots ,x_m \, \mid \, x_i= w_{i,j,k}^{-1}
x_jw_{i,j,k} , \, \, w_{i,j,k}\in W\, \rangle ,
\end{equation}
where $W=\{ w_{i,j,k}\in \mathbb F_m\, \mid \, 1\leq i,j\leq m,\,
\, 1\leq k\leq h(i,j)\}$ is a subset of elements of the free group
$\mathbb F_m$ (it is possible that
$w_{i_1,j_1,k_1}=w_{i_2,j_2,k_2}$ for $(i_1,j_1,k_1)\neq
(i_2,j_2,k_2)$), generated by free generators $x_1,\dots,x_m$ and
$h:\{1,\dots , m\}^2 \to \mathbb Z$ is some function. Such a
presentation is called a $C$-{\it presentation} ($C$, since all
relations are conjugations). Let $\varphi_W: \mathbb F_m\to G_W$
be the canonical epimorphism. The elements $\varphi_{W}(x_i)\in
G$, $1\leq i\leq m$, and the elements conjugated to them are
called the {\it $C$-generators} of the $C$-group $G$. Let
$f:G_1\to G_2$ be a homomorphism of $C$-groups. It is called a
{\it $C$-homomorphism} if the images of the $C$-generators of
$G_1$ under $f$ are $C$-generators of the $C$-group $G_2$.
$C$-groups will be considered up to $C$-isomorphisms.

Any $C$-group $G$ can be realized as the fundamental group
$\pi_1(S^4\setminus S)$ of the complement of a closed oriented
surface $S$ in the 4-dimensional sphere $S^4$ (see, for example,
\cite{Ku1}\footnote{Some other properties of $C$-groups can be
found in \cite{Ku3}, \cite{Ku}, \cite{Kuz}.}).

A $C$-presentation (\ref{zero}) is called a {\it Hurwitz
$C$-presentation of degree} $m$ if for each $i=1,\dots,m$ the word
$w_{i,i,1}$ coincides with the product $x_1\dots x_m$, and a
$C$-group $G$ is called a {\it Hurwitz $C$-group of degree} $m$ if
it possesses a Hurwitz $C$-presentation of degree $m$. In other
words, a $C$-group $G$ is a Hurwitz $C$-group of degree $m$ if
there are $C$-generators $x_1,\dots, x_m$ generating $G$ such that
the product $x_1\dots x_m$ belongs to the center of $G$. Note that
the degree of a Hurwitz $C$-group $G$ is not defined canonically
and depends on the Hurwitz $C$-presentation of $G$. Denote by
$\mathcal H$ the class of all  Hurwitz $C$-groups.

Let $\bar H$ be a Hurwitz curve of degree $m$. A Zariski -- van
Kampen presentation of $\pi_1=\pi_1(\mathbb C^2\setminus H)$
(where $\mathbb C^2=\mathbb P^2\setminus L_{\infty}$ and
$L_{\infty}$, a fibre of $\pr$, is in general position with
respect to $\bar H$) defines on $\pi_1$ a structure of a Hurwitz
$C$-group of degree $m$ (see \cite{Ku}), and in \cite{Ku}, it was
proved that any Hurwitz $C$-group $G$ of degree $m$ can be
realized as the fundamental group $\pi_1(\mathbb C^2\setminus H)$
for some Hurwitz curve $\bar H$ with singularities of the form
$w^m-z^m=0$, $\deg \bar H=2^nm$, where $n$ depends on the Hurwitz
$C$-presentation of $G$. Since we consider $C$-groups up to
$C$-isomorphisms, the class $\mathcal H$ coincides with the class
$\{ \, \pi_1(\mathbb C^2\setminus H)\, \}$ of fundamental groups
of the complements of affine Hurwitz curves.

It is easy to show that\footnote{For a group $G$ we use the
standard notation $G'$ for its commutator subgroup and $[g_1,g_2]$
for the commutator of elements $g_1,g_2\in G$.} $G/G'$ is a
finitely generated free abelian group for any $C$-group $G$. A
$C$-group $G$ is called {\it irreducible} if $G/G'\simeq \mathbb
Z$ and we say that $G$ consists of $n$ {\it irreducible
components} if $G/G'\simeq \mathbb Z^n$. If a Hurwitz $C$-group
$G$ is realized as the fundamental group $\pi_1(\mathbb
C^2\setminus H)$ of the complement of some affine Hurwitz curve
$H$, then the number of irreducible components of $G$ is equal to
the number of irreducible components of $H$. Similarly, if a
$C$-group $G$ consisting of $n$ irreducible components is realized
as $G=\pi_1(S^4\setminus S)$, then the number of connected
components of the surface $S$ is equal to $n$. \\

A free group $\mathbb F_n$ with fixed free generators is a
$C$-group and for any $C$-group $G$ the canonical $C$-epimorphism
$\nu : G\to \mathbb F_1$, sending the $C$-generators of $G$ to the
$C$-generator of $\mathbb F_1$, is well defined. Denote by $N$ its
kernel. Note that if $G$ is an irreducible $C$-group, then $N$
coincides with  $G'$.

Let $G$ be a $C$-group. The $C$-epimorphism $\nu$ induces the
following exact sequence of groups
$$1\to N/N'\to G/N'\buildrel{\nu_*}\over\longrightarrow  \mathbb F_1\to 1.$$
The $C$-generator of $\mathbb F_1$ acts on $N/N'$ by conjugation
$\widetilde x^{-1} g\widetilde{x}$, where $g\in N$ and
$\widetilde{x}$ is one of the $C$-generators of $G$. Denote by $h$
this action and by $h_{\mathbb C}$ the induced action on
$(N/N')_{\mathbb C}=(N/N')\otimes \mathbb C$. The characteristic
polynomial $\Delta (t)=\det (h_{\mathbb C}-t\text{Id})$ is called
the {\it Alexander polynomial} of the $C$-group $G$ (if the vector
space $(N/N')_{\mathbb C}$ over $\mathbb C$ is infinite
dimensional, then, by definition, the Alexander polynomial $\Delta
(t)\equiv 0$). For a Hurwitz curve $\bar H$ the Alexander
polynomial $\Delta (t)$ of the group $\pi_1=\pi_1(\mathbb
C^2\setminus H)$ is called the {\it Alexander polynomial} of $\bar
H$. Note that the Alexander polynomial $\Delta (t)$ of a Hurwitz
curve $\bar H$ does not depend on the choice of the generic line
$L_{\infty}$.

The homomorphism $\nu :\pi_1\to \mathbb F_1$ defines an infinite
unramified cyclic covering $f_{\infty}:X_{\infty}\to \mathbb
C^2\setminus H$. We have $H_1(X_{\infty},\mathbb Z)=N/N^{\prime}$
and the action of $h$ on $H_1(X_{\infty},\mathbb Z)$ coincides
with the action of a generator of the covering transformation
group of the covering $f_{\infty}$. As it follows from
\cite{Ku.O}, the\ group\ $H_1(X_{\infty},\mathbb Z)$ is finitely
generated. \\

For any $n\in \mathbb N$ denote\ by\! $\text{mod}_{n}:\mathbb
F_1\to \mu_n=\mathbb F_1/\{h^n\}$ the\ na\-tu\-ral epimorphism to
the cyclic group $\mu_n$ of degree $n$. The covering $f_{\infty}$
can be factorized through the cyclic covering $f_n:X_n\to \mathbb
C^2\setminus H$ associated with the epimorphism\
$\text{mod}_n\circ \nu$, $f_{\infty}=g_n\circ f_n$. Since a
Hur\-witz curve $\bar H$ has only analytic singularities, the
covering $f_n$ can be extended to a smooth map ${\overline
f}_n:{\overline X}_n\to \mathbb C\mathbb P^2$ branched along $\bar
H$ and, maybe, along $L_{\infty}$ (if $n$ is not a divisor of
$\deg \bar H$, then ${\overline f}_n$ is branched along
$L_{\infty}$), where $\overline X_n$ is a smooth 4-fold. The
action $h$ induces an action $\overline h_n$ on $\overline X_{n}$
and an action $\overline h_{n*}$ on $H_1(\overline X_{n},\mathbb
Z)$.

In \cite{Gr-Ku}, it was shown that any such $\overline X_n$ can be
embedded as a symplectic submanifold to a projective rational
3-fold on which the symplectic structure is given by an integer
K\"{a}hler form, and it was proved that the first Betti number
$b_1(\overline X_n)=\dim_{\mathbb C} H_1(\overline X_n,\mathbb C)$
of $\overline X_n$ is equal to $r_{n,\neq 1}$, where $r_{n,\neq
1}$ is the number of roots of the Alexander polynomial $\Delta
(t)$ of the curve $\bar H$ which are $n$-th roots of unity not
equal to 1. \\

In \cite{Gr-Ku}, properties of the Alexander polynomials of
Hurwitz curves were investigated. In particular, it was proved
that if $\bar H$ is a Hurwitz curve of degree $d$ consisting of
$n$ irreducible components and $\Delta (t)$ is its Alexander
polynomial, then
\begin{itemize}
\item[($i$)] $\Delta (t)\in \mathbb Z[t]$; \item[($ii$)]
$\Delta(0)=\pm 1$; \item[($iii$)] the roots of $\Delta (t)$  are
$d$-th roots of unity; \item[($iv$)] the action of $h_{\mathbb C}$
on $(N/N')\otimes \mathbb C$ is semisimple;
\item[($v$)] $\Delta (t)$ is a divisor of the polynomial
$(t-1)(t^d-1)^{d-2}$;
\item[($vi$)] the multiplicity of
the root $t=1$ of $\Delta (t)$ is equal to $n-1$;
\item[($vii$)] if $n=1$, then $\Delta (1)=1$ and $\deg \Delta (t)$
is an even number.
\end{itemize}

The main results of the article are the following theorems
\footnote{Recall that the class of Hurwitz $C$-groups coincides
with the class of fundamental groups of the complements of affine
Hurwitz curves. Therefore to speak about the Alexander polynomials
of Hurwitz curves is the same as to speak about the Alexander
polynomials of Hurwitz $C$-groups. Hence the results of the
article are formulated in terms of Hurwitz $C$-groups, since their
proves are purely algebraic.}.

\begin{thm} \label{irred} A polynomial $P(t)\in \mathbb Z[t]$ is the Alexander polynomial of
an irreducible Hurwitz $C$-group $G$ iff the roots of $P(t)$  are
roots of unity and $P(1)=1$.
\end{thm}

\begin{thm} \label{nirred} Let a polynomial $P(t)=(-1)^{m} t^m+\sum_{i=0}^{m-1}a_it^i\in \mathbb Z[t]$ has
the following properties:
\begin{itemize}
\item[($i$)] the roots of $P(t)$  are roots of unity;
 \item[($ii$)] if $\zeta$ is a primitive $p^k$-th root of unity, where $p$ is a prime number, then
 the multiplicity of the root $t=\zeta$ of the polynomial $P(t)$ is not greater than the multiplicity
 of its root $t=1$.
\end{itemize}
Then $P(t)$ is the Alexander polynomial of some Hurwitz $C$-group
$G$.
\end{thm}

\begin{thm} \label{two} The polynomial $P(t)=(-1)^{n+k}(t-1)^n(t+1)^k$
is the Alexander polynomial of some Hurwitz $C$-group $G$ iff
$n\geq k$.
\end{thm}

The proof of Theorem \ref{irred} is given in section 1. Section 2
is devoted to the proof of Theorem \ref{nirred}. In spite of the
proves of these theorems use almost the same ideas, the proves are
given independently in order to have more clear exposition. The
proof of Theorem \ref{two} is given in section 3. These theorems
allow to formulate the following
\begin{conj} \label{nirreduc} A polynomial
$P(t)=(-1)^{m} t^m+\sum_{i=0}^{m-1}a_it^i\in \mathbb Z[t]$ is the
Alexander polynomial of some Hurwitz $C$-group $G$ iff it
satisfies the following conditions:
\begin{itemize}
\item[($i$)] the roots of $P(t)$  are
 roots of unity;
 \item[($ii$)] if $\zeta$ is a primitive $p^k$-th root of unity,
 where $p$ is a prime number, then
 the multiplicity of the root $t=\zeta$ of the polynomial $P(t)$
 is not greater than the multiplicity of its root $t=1$.
\end{itemize}
\end{conj}

In section 4, to compare the case of Hurwitz $C$-groups with the
general case of $C$-groups, we give an example (see Example
\ref{ex2}) of a $C$-group consisting of two irreducible components
whose Alexander polynomial $\Delta(t)=(1-t)(t+1)^2$ (compare it
with Theorem \ref{two}), and  an example (see Example \ref{ex1})
of a $C$-group consisting also of two irreducible components whose
Alexander polynomial $\Delta(t)=(t-1)^2$ (compare it with property
($vi$) mentioned above, see also Theorem 5.9 in \cite{Gr-Ku}).

\section{Alexander polynomials of irreducible Hurwitz $C$-groups}

In this section, we prove Theorem \ref{irred}.

In \cite{Gr-Ku}, it was proved that if a polynomial $P(t)\in
\mathbb Z[t]$ is the Alexander polynomial of an irreducible
Hurwitz $C$-group $G$, then the roots of $P(t)$  are roots of
unity and $P(1)=1$. Therefore to prove Theorem \ref{irred}, it
suffices to prove the inverse statement.

Consider a polynomial $P(t)\in \mathbb Z[t]$ such that the roots
of $P(t)$  are roots of unity and $P(1)=1$.

In the beginning, let us show that it suffices to prove that for
any polynomial $\Psi (t)$ such that
\begin{itemize}
\item[($i$)] the roots of $\Psi (t)$ are roots of unity,
\item[($ii$)] $\Psi (t)$ has not
multiple roots,
 \item[($iii$)] $\Psi (1)=1$, \end{itemize} there
exists an irreducible Hurwitz $C$-group with the Alexander
polynomial $\Delta(t)=\Psi (t)$. Indeed, each polynomial $P(t)\in
\mathbb Z[t]$ can be factorized into the product
$P(t)=\prod_{i}\Psi_i(t)$, where each factor $\Psi_i(t)\in \mathbb
Z[t]$ has not multiple roots. Since the roots of $P(t)$ are roots
of unity, $P(1)=1$, and $\Psi_i(1)\in \mathbb Z$, we have
$\Psi_i(1)=1$ for each $i$. Next, if $\Delta_1(t)$ and
$\Delta_2(t)$ are the Alexander polynomials of two irreducible
Hurwitz $C$-groups $G_1$ and $G_2$ given by Hurwitz
$C$-presentations $G_i=\langle x_{1,i},\dots,x_{m_i,i}\, \mid \,
\mathcal R_i \rangle $ of degrees $m_i$, $i=1,2$, then, by
Proposition 5.12 of \cite{Gr-Ku}, the Alexander polynomial of a
Hurwitz product $G=G_1\diamond G_2$ is equal to
$\Delta_1(t)\Delta_2(t)$, where the Hurwitz product $G$ of  $G_1$
and $G_2$ is an irreducible Hurwitz $C$-group given by
$C$-presentation
\begin{align*}  G= & \langle x_{1,i}\dots ,x_{m_i,i}, i=1,2 \mid \mathcal R_1\cup
\mathcal R_2,\, x_{m_1,1}=x_{m_2,2}, \\ &
[x_{j,i},(x_{1,{\overline i}}\dots x_{m_{\overline i},{\overline
i}})^{m_i}]=1 \, \text{for}\, j=1,\dots,m_i-1, i=1,2\rangle
,\end{align*} where $\overline i=\{ 1,2\}\setminus \{ i\}$ (recall
that it follows from the set of relations $\mathcal R_{\overline
i}$ that the elements $x_{1,\overline i},\dots, x_{m_{\overline
i},\overline i}$ commute with the product $x_{1,{\overline
i}}\dots x_{m_{\overline i},{\overline i}}$).

Fix one of the polynomials $\Psi_i(t)=\Psi (t)$. Denote by $k$ the
smallest positive integer such that all roots of $\Psi (t)$ are
$k$-th roots of unity. Since $\Psi (t)$ has not multiple roots,
$\Psi (t)$ is a divisor of the polynomial $t^k-1$.

To prove the existence of an irreducible Hurwitz $C$-group whose
Alexander polynomial coincides with $\Psi (t)$, consider the
quotient ring $M=M_{\Psi}=\mathbb Z[t]/(\Psi (t))$ of the ring
$\mathbb Z[t]$.

The ring $M$ (considered as an abelian group with respect to
additive operation) is a free abelian group freely generated by
the elements $t_0=t^0,\dots,t_{d-1}=t^{d-1}$, where $d=\deg \Psi
(t)$. Let $h_0\in \text{End}_{\mathbb Z}M$ given by
$h_0(Q(t))=tQ(t)$ and $i:\mathbb F_1\to \text{End}_{\mathbb Z}M$
given by $i(x_0)=h_0$, where $x_0$ is a generator of the free
group $\mathbb F_1$. Obviously, $h_0\in \text{Aut}\, M$. It is not
hard to show (see, for example, \cite{L}, Chapter XV) that
$\det(h_{0,\mathbb C}-t\text{Id})=\Psi (t)$, where $h_{0,\mathbb
C}\in \text{Aut}(M\otimes \mathbb C)$ is induced by $h_0$.

Since $h_0\in \text{Aut}\, M$, we can consider (below, we use the
multiplicative notation for the group operation in $M$) the
semi-direct product
$$
\begin{array}{ll}
G_{\Psi}=M\leftthreetimes\mathbb F_1\simeq &  \\
 \langle
t_0,\dots,t_{d-1},x_0\mid & [t_i,t_j]=1\, \qquad \qquad
\text{for}\,\, 0\leq i,j\leq
d-1,\\ & x_0^{-1}t_ix_0=t_{i+1},\,\, \qquad \text{for}\,\, i=0,\dots, d-2,\\
& 
x_0^{-1}t_{d-1}x_0=\displaystyle
\prod_{i=0}^{d-1}t_i^{-a_i}\rangle ,
\end{array}
$$
where $a_i$ are taken from the equality $\Psi
(t)=t^d+\sum_{i=0}^{d-1}a_it^i$.

Let $\nu : G_{\Psi}\to \mathbb F_1$ be the canonical epimorphism,
$\ker \nu=M\subset G_{\Psi}$.  Obviously,  $h_0(g)=x_0^{-1}gx_0$
for $g\in M$.

For $i\geq d$ denote by $t_i=x_0^{-i}t_0x^i$. Evidently, $x_0^k$
belongs to the center of $G_{\Psi}$, since $\Psi (t)$ is a divisor
of the polynomial $t^k-1$ and therefore
$$t^{k+i}\equiv t^i \mod \Psi (t)$$ for all $i\geq 0$.
Thus, we have $t_i=t_{k+i}$ for all $i\geq 0$.

For $P(t)=\sum_{i=0}^nc_it^i\in \mathbb Z[t]$ denote by
$x_{P(t)}=x_0\prod_{i=0}^nt_i^{c_i}$. Obviously,
\begin{equation} \label{equ1} x_{P(t)}=x_{Q(t)}\qquad \text{if}\,\, P(t)\equiv Q(t)\, \text{mod}\,
\Psi (t)
\end{equation} and it is easy to check that
\begin{equation} \label{equ2}
x_0^{-1}x_{P(t)}x_0=x_{tP(t)},
\end{equation}
and
\begin{equation} \label{equ3}
(\prod_{i=0}^nt_i^{c_i})x_{Q(t)}(\prod_{i=0}^nt_i^{c_i})^{-1}=x_{(t-1)P(t)+Q(t)}
\end{equation}
for $P(t)=\sum_{i=0}^nc_it^i$ and for any $Q(t)$. In particular,
\begin{equation} \label{rel1}
x_{t^{i+1}}=x_0^{-1}x_{t^i}x_0
\end{equation}
for all $i\geq 0$.

Since $\Psi (1)=1$, one can find a polynomial $P(t)\in \mathbb
Z[t]$ such that $\Psi (t)=(t-1)P(t)+1$. Therefore it follows from
(\ref{equ1}) and (\ref{equ3}) that there is $g_{1,0}\in G_{\Psi}$
such that $x_{t^0}=g_{1,0}^{-1}x_0g_{1,0}$. Consequently, by
(\ref{equ2}), the elements $x_0,x_{t^0},\dots,x_{t^{d-1}}$ are
conjugated to each other and therefore $x_{t^i}^k$ belongs to the
center of $G_{\Psi}$ for each $i$.

Since $t_i=x_0^{-1}x_{t^i}$ and $x_0,t_0,\dots, t_{d-1}$ generate
$G_{\Phi_k}$, the elements $x_0,x_{t^0},\dots, x_{t^{d-1}}$ also
generate $G_{\Psi}$. Let $w_{1,0}(x_0,x_{t^0},\dots,x_{t^{d-1}})$
be a word in letters $x_0,x_{t^0},\dots,x_{t^{d-1}}$ and their
inverses representing the element $g_{1,0}$. Consider a $C$-group
$\widetilde G$ given by $C$-presentation
\begin{equation} \label{group}
\begin{array}{ll}
\widetilde G=
 \langle
x_1,\dots,x_{d+1}\, \, \mid & x_{i+1}=x_1^{-1}x_ix_1\,  \qquad
\text{for}\,\, 2\leq i\leq d,\\ & x_1^{-k}x_ix_1^k=x_i \quad
\qquad \text{for}\,\, 2\leq i\leq d+1, \\ &
x_2=w_{1,0}^{-1}(x_1,\dots,x_{d+1})x_1w_{1,0}(x_1,\dots,x_{d+1})\,
\rangle .
\end{array}
\end{equation}

First of all, note that $\widetilde G$ is an irreducible
$C$-group, since all $C$-generators of $\widetilde G$ are
conjugated to each other. Next, the element $x_1^k$ belongs to the
center of $\widetilde G$. Therefore $x_i^k=x_1^k$ for $2\leq i\leq
d+1$, since the elements  $x_i^k$ are also conjugated to $x_1^k$
in $\widetilde G$. Hence the product $x_1^k\dots x_{d+1}^k$
belongs to the center of $\widetilde G$. Therefore, by Lemma 5.11
of \cite{Gr-Ku}, $\widetilde G$ is a Hurwitz $C$-group of degree
$k(d+1)$. (Indeed, one can consider another $C$-presentation of
$\widetilde G$ with $C$-generators $\overline x_{i,j}$, $1\leq
i\leq d+1$, $1\leq j\leq k$, satisfying the following defining
relations:
$$ \begin{array}{rcll} \overline
x_{i,j_1}& = & \overline x_{i,j_2} &  \text{for all}\, \,  i, j_1,
j_2; \\
x_{1,1}^{-1}x_{i,1}x_{1,1} & = & x_{i+1,1} & \text{for}\, \, 2\leq
i\leq d; \\
x_{1,1}^{-k}x_{i,1}x_{1,1}^k & = & x_{i,1} & \text{for}\,
\, 2\leq i\leq d+1; \\
\displaystyle [ \prod_{i=1}^{d+1}\prod_{j=1}^k\overline x_{i,j},
\overline x_{i,j}] & = & 1 & \text{for all}\, \, i\, \,
\text{and}\, \, j;
\end{array} $$
$$x_{2,1}=w_{1,0}^{-1}(x_{1,1},\dots
,x_{d+1,1})x_{1,1}w_{1,0}(x_{1,1},\dots ,x_{d+1,1}).$$ Obviously,
this presentation becomes a Hurwitz $C$-presentation after
renumbering the generators).

Denote by $\widetilde \nu: \widetilde G\to \mathbb F_1$ the
canonical $C$-epimorphism and by $\widetilde h_1$ the action of
the $C$-generator of $\mathbb F_1$ on $\widetilde N=\ker
\widetilde \nu$ given by $\widetilde h_1(g)=x_1^{-1}gx_1$ for
$g\in \widetilde N$.

It is easy to see that a map $\widetilde f$, sending the set $\{
x_i\mid i=1,\dots, d+1\}$ of generators of $\widetilde G$ to the
set $\{ x_0, x_{t^i}\mid i=0,\dots, d-1\}$ of generators of
$G_{\Psi}$ as follows: $\widetilde f(x_1)=x_0$ and $\widetilde
f(x_i)=x_{t^{i-2}}$ for $2\leq i\leq d+1$, can be extended to an
epimorphism $\widetilde f:\widetilde G\to G_{\Psi}$. Obviously, we
have $\widetilde f(\widetilde N)=M$ and $h_0(\widetilde
f(g))=\widetilde f(\widetilde h_1(g))$ for all $g\in \widetilde
N$. Consider the restriction $\widetilde f_{\mid \widetilde
N}:\widetilde N\to M$ of $\widetilde f$ to $\widetilde N$. Since
$M$ is an abelian group, the epimorphism $\widetilde f_{\mid
\widetilde N}$ can be decomposed into the composition $\widetilde
f_{\mid \widetilde N}=\widetilde f'\circ \widetilde j$, where
$\widetilde j:\widetilde N \to \widetilde N/\widetilde N'$ is the
canonical epimorphism and $\widetilde f':\widetilde N/\widetilde
N'\to M$ is the epimorphism induced by $\widetilde f_{\mid
\widetilde N}$. Denote by $\widetilde K=\ker \widetilde f'$.

It follows from \cite{Ku.O} that the abelian group $\widetilde
N/\widetilde N'$ is finitely generated.  Therefore the group
$\widetilde K$ is also a finitely generated abelian group. Let
$g_1,\dots,g_n$ be a set of elements of $\widetilde N$ such that
their images $\widetilde j(g_1)$, $\dots$, $\widetilde j(g_n)$ in
$\widetilde N/\widetilde N'$ generate the group $\widetilde K$ and
let a word $w_i(x_1,\dots,x_{d+1})$, $1\leq i\leq n$, in the
letters
 $x_1^{\pm 1},\dots,x_{d+1}^{\pm
1} $ represent the element $g_i$. Consider a $C$-group $G$ given
by $C$-presentation
\begin{equation} \label{group2}
\begin{array}{ll}
G=
 \langle
x_1,\dots,x_{d+1} \mid \mathcal R, &
x_j=w_{i}^{-1}(x_1,\dots,x_{d+1})x_jw_{i}(x_1,\dots,x_{d+1}) \\
& 1\leq j\leq d+1,\, 1\leq i\leq n\,  \rangle ,
\end{array}
\end{equation}
where $\mathcal R$ is the set of defining relations in
presentation (\ref{group}). The group $G$ is also a Hurwitz
$C$-group. Moreover, it follows from the defining relations for
$G$ that the epimorphism $\widetilde f$ can be decomposed into the
composition $\widetilde f=f\circ r$ where the $C$-epimorphism $r$
sends the generators $x_i$, $1\leq i\leq d+1$ of the group
$\widetilde G$ to the generators $x_i$, $1\leq i\leq d+1$ of the
group $G$ and $f:G\to G_{\Psi}$ is an epimorphism such that
$f(x_1)=x_0$ and $f(x_i)=x_{t^{i-2}}$ for $2\leq i\leq d+1$. It
follows from (\ref{group2}) that the elements $r(g_i)$, $1\leq
i\leq n$, belong to the center of the group $G$.

Denote by $\nu: G\to \mathbb F_1$ the canonical $C$-epimorphism
and by $\widetilde h$ the action of the $C$-generator of $\mathbb
F_1$ on $N=\ker \nu$ given by $\widetilde h(g)=x_1^{-1}gx_1$ for
$g\in N$.  Obviously, we have $f(N)=M$ and $
h_0(f(g))=f(\widetilde h(g))$ for all $g\in N$.

The restriction $f_{\mid N}: N\to M$ of $f$ to $N$ can be also
decomposed into the composition $f_{\mid  N}=f'\circ j$, where $
j: N \to N/N'$ is the canonical epimorphism and $f': N/N'\to M$ is
the epimorphism induced by $f_{\mid  N}$. Denote by $K=\ker f'$.
Obviously, we have $r'(\widetilde K)=K$, where $r':\widetilde
N/\widetilde N'\to N/N'$ is the epimorphism induced by $r$. Let
$h\in \text{Aut}\, N/N'$ be the automorphism induced by
$\widetilde h$.

We have the following commutative diagram.
\\

\begin{picture}(0,0)(0,0)
\put(80,-5){$\widetilde N$} \put(95,-0){\vector(1,0){45}}
\put(115,5){$\widetilde j$} \put(148,-5){$\widetilde N/\widetilde
N'$} \put(80,-90){$N$} \put(95,-85){\vector(1,0){45}}
\put(115,-97){$j$} \put(146,-90){$N/N'$}
\put(85,-10){\vector(0,-1){65}} \put(78,-45){$r$}
\put(158,-10){\vector(0,-1){65}} \put(148,-45){$r'$}
\put(120,-30){$\widetilde K$} \put(120,-65){$K$}
\put(130,-75){$\searrow$} \put(132,-15){$\nearrow$}
\put(123,-33){\vector(0,-1){18}} \put(114,-45){$r'$}
\put(230,-45){$M$} \put(180,-7){\vector(3,-2){45}}
\put(200,-17){$\widetilde f'$} \put(200,-77){$f'$}
\put(180,-80){\vector(3,2){45}}  \put(330,-45){$(*)$}
\end{picture}
\\
\\
\\
\\
\\
\\
\\
\\

The elements $\widetilde j(g_1),\dots,\widetilde j(g_n)$ generate
the group $\widetilde K$. Since the homomorphisms $\widetilde j,
j, r$, and $r'$ in commutative diagram ($*$) are epimorphisms, the
group $K$ is generated by the elements
$j(r(g_1)),\dots,j(r(g_n))$.

By construction of the group $G$, the elements
$r(g_1),\dots,r(g_n)$ belong to the center of $G$. Therefore
$h(j(r(g_i)))=j(r(g_i))$, that is, $h_{\mid K}=\text{Id}$. But
$t=1$ is not a root of the Alexander polynomial $\Delta
(t)=\det(h_{\mathbb C}-t\text{Id})$ of $G$, since $G$ is an
irreducible $C$-group. Therefore $K$ belongs to the kernel of the
natural homomorphism $c:N/N'\to (N/N')_{\mathbb C}=(N/N')\otimes
\mathbb C$ which coincides with the subgroup $\text{Tors}\, N/N'$
of $N/N'$ consisting of the elements of finite orders. On the
other hand, $\text{Tors}\, N/N'\subset \ker f'$, since $M$ is a
free abelian group. Therefore $K=\text{Tors}\, N/N'$ and, since
$f'$ is an epimorphism, $f'_{\mathbb C}:(N/N')_{\mathbb C}\to
M_{\mathbb C}$ is an isomorphism. Thus, $\Delta (t)=\Psi (t)$,
since $h_0(f'(g))=f'(h(g))$ for all $g\in N/N'$. \qed

\section {The Alexander polynomials of Hurwitz $C$-groups
consisting of several irreducible components} In this section,
Theorem \ref{nirred} will be proved.

Consider a polynomial $P(t)=(-1)^{n+\deg P_0}(t-1)^nP_0(t) \in
\mathbb Z[t]$ such that $P_0(1)\neq 0$, all roots of $P_0(t)$ are
roots of unity, and if $\zeta$ is a primitive $p^k$-th root of
unity, where $p$ is a prime number, then the multiplicity of the
root $t=\zeta$ of the polynomial $P_0(t)$ is not greater than $n$.
The polynomial $P_0(t)$ can be factorized into the product of
cyclotomic polynomials. It is well known (see, for example, Lemma
5.3 of \cite{Gr-Ku}) that for a cyclotomic polynomial $\Phi_k(t)$
the value $\Phi_k(1)\neq 1$ iff $k=p^m$ for some prime number $p$.
Therefore one can find a factorization
$P(t)=\prod_{i}^{n+1}\Psi_i(t)$, where for $i=1,\dots, n$ each
factor $\Psi_i(t)=(-1)^{\deg \Psi_{i,0}}(1-t)\Psi_{i,0}(t)\in
\mathbb Z[t]$ has not multiple roots and $\Psi_{n+1}(t)\in \mathbb
Z[t]$ is a polynomial such that $\Psi_{n+1}(1)=1$.

By Theorem \ref{irred}, there is a Hurwitz $C$-group whose
Alexander polynomial coincides with $\Psi_{n+1}(t)$. Besides, by
Proposition 5.12 of \cite{Gr-Ku}, the Alexander polynomial of a
Hurwitz product $G=G_1\diamond G_2$ of two Hurwitz $C$-groups is
equal to the product of their Alexander polynomials. Therefore to
prove Theorem \ref{nirred}, it suffices to prove the existence of
a Hurwitz $C$-group $G$ whose Alexander polynomial $\Delta
(t)=\Psi (t)$, where $\Psi (t) =(-1)^{\deg \Psi}(1-t)\Psi (t)\in
\mathbb Z[t]$ has no multiple roots and all its roots are roots of
unity.

Fix such a polynomial $\Psi (t)$. Let $d=\deg \Psi (t)$. Denote by
$k$ the smallest positive integer such that all roots of $\Psi
(t)$ are $k$-th roots of unity. Since $\Psi (t)$ has not multiple
roots, $\Psi (t)$ is a divisor of the polynomial $t^k-1$.
Therefore, $\Psi (0)=\pm 1$ and $\Psi (t)=\pm t^d+\dots$, that is,
its leading coefficient is equal to $\pm 1$.

As in the proof of Theorem \ref{irred}, to prove the existence of
a Hurwitz $C$-group whose Alexander polynomial coincides with
$\Psi (t)$, consider the quotient ring $M=M_{\Psi}=\mathbb
Z[t]/(\Psi (t))$. Since the leading coefficient of $\Psi (t)$ is
equal to $\pm 1$, the ring $M$ (considered as an abelian group
with respect to additive operation) is a free abelian group freely
generated by the elements $t_0=t^0,\dots,t_{d-1}=t^{d-1}$. As
above, let $h_0\in \text{End}_{\mathbb Z}M$ given by
$h_0(Q(t))=tQ(t)$ and $i:\mathbb F_1\to \text{End}_{\mathbb Z}M$
given by $i(x_0)=h_0$, where $x_0$ is a generator of the free
group $\mathbb F_1$. We have $h_0\in \text{Aut}\, M$, since $\Psi
(0)=\pm 1$ and $\det(h_{0,\mathbb C}-t\text{Id})=\Psi (t)$.

Since $h_0\in \text{Aut}\, M$, we can consider the semi-direct
product
$$
\begin{array}{ll}
G_{\Psi}=M\leftthreetimes\mathbb F_1\simeq &  \\
 \langle
t_0,\dots,t_{d-1},x_0\mid & [t_i,t_j]=1\, \qquad \qquad
\text{for}\,\, 0\leq i,j\leq
d-1,\\ & x_0^{-1}t_ix_0=t_{i+1},\,\, \qquad \text{for}\,\, i=0,\dots, d-2,\\
& 
x_0^{-1}t_{d-1}x_0=\displaystyle
\prod_{i=0}^{d-1}t_i^{-a_i}\rangle ,
\end{array}
$$
where $a_i$ are coefficients of the polynomial $\Psi
(t)=t^d+\sum_{i=0}^{d-1}a_it^i$.

Let $\nu_0 : G_{\Psi}\to \mathbb F_1$ be the canonical
epimorphism, $\ker \nu_0=M\subset G_{\Psi}$. We have
$h(g)=x_0^{-1}gx_0$ for all $g\in M$.

For $i\geq d$ denote by $t_i=x_0^{-i}t_0x^i$. Obviously, $x_0^k$
belongs to the center of $G_{\Psi}$, since $\Psi (t)$ is a divisor
of the polynomial $t^k-1$ and therefore
$$t^{k+i}\equiv t^i \mod \Psi (t)$$ for all $i\geq 0$.
Thus, we have $t_i=t_{k+i}$ for all $i\geq 0$.

Denote by $x_{t^i}=x_0t_i$. We have
\begin{equation} \label{equ11}
x_0^{-1}x_{t^i}x_0=x_{t^{i+1}}.
\end{equation}

Since each $t_i$ belongs to the normal abelian subgroup $M$ of
$G_{\Psi}$, we have the equality
$x_0^{-1}gx_0=x_{t^i}^{-1}gx_{t^i}$ for all $g\in M$. Therefore
all the elements $x_{t^i}^k$ also belong to the center of
$G_{\Psi}$, since $x_0^k$ belongs to its center.

We have $t_i=x_0^{-1}x_{t^i}$ and the elements $x_0,t_0,\dots,
t_{d-1}$ generate $G_{\Phi}$. Therefore the elements
$x_0,x_{t^0},\dots, x_{t^{d-1}}$ also generate $G_{\Psi}$.

Consider a $C$-group $\widetilde G$ given by $C$-presentation
\begin{equation} \label{group3}
\begin{array}{ll}
\widetilde G=
 \langle
x_1,\dots,x_{k+1}\, \, \mid &
x_{i+1}=x_1^{-1}x_ix_1\, \, \text{for}\,\, 2\leq i\leq k,\, \,  x_1^{-1}x_{k+1}x_1=x_2, \\
& [x_i^k,x_j]=1  \quad \text{for}\,\, 2\leq i\leq k,\, \, 1\leq
j\leq k+1 \rangle .
\end{array}
\end{equation}

First of all, notice that $\widetilde G$ is a $C$-group consisting
of two irreducible components. Next, all the elements $x_i^k$,
$i=1,\dots,k+1$, belong to the center of $\widetilde G$. Thus, the
product $x_1^k\dots x_{k+1}^k$ belongs to the center of
$\widetilde G$. Therefore, as in the proof of Theorem \ref{irred},
one can show that $\widetilde G$ is a Hurwitz $C$-group of degree
$k(k+1)$.

Denote, as above, by $\widetilde \nu: \widetilde G\to \mathbb F_1$
the canonical $C$-epimorphism and by $\widetilde h_1$ the action
of the $C$-generator of $\mathbb F_1$ on $\widetilde N=\ker
\widetilde \nu$ given by $\widetilde h_1(g)=x_1^{-1}gx_1$ for all
$g\in \widetilde N$.

It is easy to see that a map $\widetilde f$, sending the set $\{
x_i\mid i=1,\dots, k+1\}$ of generators of $\widetilde G$ to the
set $\{ x_0, x_{t^i}\mid i=0,\dots, k-1\}$ of generators of
$G_{\Psi}$ as follows: $\widetilde f(x_1)=x_0$ and $\widetilde
f(x_i)=x_{t^{i-2}}$ for $2\leq i\leq k+1$, can be extended to an
epimorphism $\widetilde f:\widetilde G\to G_{\Psi}$. Obviously, we
have $\widetilde f(\widetilde N)=M$ and $h_0(\widetilde
f(g))=\widetilde f(\widetilde h_1(g))$ for all $g\in \widetilde
N$.

As in the proof of Theorem \ref{irred}, the restriction
$\widetilde f_{\mid \widetilde N}:\widetilde N\to M$ of
$\widetilde f$ to $\widetilde N$ can be decomposed into the
composition $\widetilde f_{\mid \widetilde N}=\widetilde f'\circ
\widetilde j$, where $\widetilde j:\widetilde N \to \widetilde
N/\widetilde N'$ is the canonical epimorphism and $\widetilde
f':\widetilde N/\widetilde N'\to M$ is the epimorphism induced by
$\widetilde f_{\mid \widetilde N}$. Denote by $\widetilde K=\ker
\widetilde f'$.

The groups $\widetilde N/\widetilde N'$ and $\widetilde K$ are
finitely generated abelian groups. Let $g_1,\dots,g_n\in
\widetilde N$ be  elements such that their images $\widetilde
j(g_1),\dots,\widetilde j(g)_n$ in $\widetilde N/\widetilde N'$
generate the group $\widetilde K$ and let a word
$w_i(x_1,\dots,x_{k+1})$, $1\leq i\leq n$, in letters
$x_1,\dots,x_{k+1}$ and their inverses represents the element
$g_i$. Consider a $C$-group $G$ given by $C$-presentation
\begin{equation} \label{group4}
\begin{array}{ll}
G=
 \langle
x_1,\dots,x_{k+1} \mid \mathcal R, &
x_j=w_{i}^{-1}(x_1,\dots,x_{k+1})x_jw_{i}(x_1,\dots,x_{k+1}) \\
& 1\leq j\leq k+1,\, 1\leq i\leq n\,  \rangle ,
\end{array}
\end{equation}
where $\mathcal R$ is the set of defining relations in
presentation (\ref{group3}). The group $G$ is also a Hurwitz
$C$-group consisting of two irreducible components. Moreover, by
construction of $G$, the epimorphism $\widetilde f$ can be
decomposed into the composition $\widetilde f=f\circ r$ where the
$C$-epimorphism $r$ sends the generators $x_i$, $1\leq i\leq d+1$,
of the group $\widetilde G$ to the generators $x_i$, $1\leq i\leq
k+1$, of the group $G$ and $f:G\to G_{\Psi}$ is an epimorphism
such that $f(x_1)=x_0$ and $f(x_i)=x_{t^{i-2}}$ for $2\leq i\leq
k+1$. It follows from (\ref{group4}) that the elements $r(g_i)$,
$1\leq i\leq n$, belong to the center of the group $G$.

Let $\nu: G\to \mathbb F_1$ be the canonical $C$-epimorphism and
let $\widetilde h$ be an action of the $C$-generator of $\mathbb
F_1$ on $N=\ker \nu$ given by $\widetilde h(g)=x_1^{-1}gx_1$ for
$g\in N$. Obviously, we have $f(N)=M$ and $ h_0(f(g))=f(\widetilde
h(g))$ for all $g\in N$.

The restriction $f_{\mid N}:N\to M$ of $f$ to $N$ also can be
decomposed into the composition $f_{\mid N}=f'\circ j$, where $ j:
N \to N/N'$ is the canonical epimorphism and $f': N/N'\to M$ is
the epimorphism induced by $f_{\mid N}$. Denote by $K=\ker f'$.
Obviously, we have $r'(\widetilde K)=K$, where $r':\widetilde
N/\widetilde N'\to N/N'$ is the epimorphism induced by $r$. Let
$h\in \text{Aut}\, N/N'$ be the automorphism induced by
$\widetilde h$.

We have the following commutative diagram.
\\

\begin{picture}(20,0)(20,0)
\put(80,-5){$\widetilde N$} \put(95,-0){\vector(1,0){45}}
\put(115,5){$\widetilde j$} \put(148,-5){$\widetilde N/\widetilde
N'$} \put(184,-0){\vector(1,0){72}} \put(214,3){$c$}
 \put(260,-5){$(\widetilde N/\widetilde N')_{\mathbb C}$}
\put(80,-90){$N$} \put(95,-85){\vector(1,0){45}}
\put(115,-97){$j$} \put(146,-90){$N/N'$}
\put(184,-85){\vector(1,0){72}} \put(214,-92){$c$}
\put(260,-90){$(N/N')_{\mathbb C}$}
\put(85,-10){\vector(0,-1){65}} \put(78,-45){$r$}
\put(158,-10){\vector(0,-1){65}} \put(148,-45){$r'$}
\put(120,-30){$\widetilde K$} \put(120,-65){$K$}
\put(130,-75){$\searrow$} \put(132,-15){$\nearrow$}
\put(123,-33){\vector(0,-1){18}} \put(114,-45){$r'$}
\put(230,-45){$M$} \put(245,-40){\vector(1,0){20}}
\put(255,-37){$c$} \put(270,-45){$M_{\mathbb C}$}
\put(340,-45){$(**)$} \put(180,-7){\vector(3,-2){45}}
\put(210,-24){$\widetilde f'$} \put(210,-70){$f'$}
\put(180,-80){\vector(3,2){45}} \put(275,-10){\vector(0,-1){20}}
\put(280,-22){$\widetilde f'_{\mathbb C}$}
\put(275,-75){\vector(0,1){22}} \put(280,-68){$f'_{\mathbb C}$}
\end{picture}
\\
\\
\\
\\
\\
\\
\\
\\

By Theorem 5.9 of \cite{Gr-Ku}, the value $t=1$ is a root of
multiplicity one of the Alexander polynomial $\Delta
(t)=\det(h_{\mathbb C}-t\text{Id})$ of the group $G$, since the
Hurwitz $C$-group $G$ consists of two irreducible component.
Besides, $t=1$ is also a root of multiplicity one of the
characteristic polynomial $\det ((h_{0,\mathbb
C}-t\text{Id})=(-1)^d\Psi (t)$. Therefore the eigen-spaces
$((N/N')_{\mathbb C})_1\subset (N/N')_{\mathbb C}$ and
$(M_{\mathbb C})_1\subset M_{\mathbb C}$ corresponding to the
eigen-value 1 are one-dimensional. Next, $f_{\mathbb C}$ is an
epimorphism such that $f'_{\mathbb C}(h_{\mathbb
C}(v))=h_{0,\mathbb C}(f'_{\mathbb C}(v))$ for any $v\in
(N/N')_{\mathbb C}$. Therefore $((N/N')_{\mathbb C})_1\not\subset
\ker f'_{\mathbb C}$.

The elements $\widetilde j(g_1),\dots,\widetilde j(g_n)$ generate
the group $\widetilde K$. Since the homomorphisms $\widetilde f,
f, r$, and $r'$ in commutative diagram ($**$) are epimorphisms,
the group $K$ is generated by the elements
$j(r(g_1)),\dots,j(r(g_n))$.

By construction of the group $G$, the elements
$r(g_1),\dots,r(g_n)$ belong to the center of $G$. Therefore,
$h(j(r(g_i)))=j(r(g_i))$, that is, $h_{\mid K}=\text{Id}$, and
consequently, $c(K)=\ker f'_{\mathbb C}\cap ((N/N')_{\mathbb C})_1
=0$. We have $K\subset \ker c=\text{Tors}\, N/N'$. On the other
hand, $\text{Tors}\, N/N'\subset K$, since $M$ is a free abelian
group. Therefore, $K=\text{Tors}\, N/N'$ and $f'_{\mathbb C}$ is
an equivariant isomorphism. Thus, $\Delta(t)= \Psi (t)$.
 \qed

\section {Hurwitz $C$-groups with the Alexander polynomials
$\Delta (t)=(-1)^{n+k}(t-1)^n(t+1)^k$}

This section is devoted to the proof of Theorem \ref{two}.

Consider a pair $A=(M,h)$, where $M$ is a finitely generated free
abelian group and $h\in \text{Aut}\, M$. In particular, set
$A_{\pm}=(\mathbb Z,\pm\text{Id})$ and  $A_{+-}=(\mathbb Z\oplus
\mathbb Z,h_{\pm})$, where, in a base $e_1,e_2$ of $\mathbb
Z\oplus \mathbb Z$, the automorphism $h_{\pm}$ is given by
$h_{\pm}(e_1)=e_2$ and $h_{\pm}(e_2)=e_1$. For two pairs $A=(M,h)$
and $A'=(M',h')$ the equality $A=A'$ will mean that $A$ and $A'$
are isomorphic, that is, there is an isomorphism $g:M\to M'$ such
that $h=g^{-1}\circ h'\circ g$. For $A$ and $A'$, put $A\oplus
A'=(M\oplus M',h\oplus h')$ and denote by $nA$ the direct sum of
$n$ copies of $A$.

\begin{prop} \label{decomposition} Let $A=(M,h)$ be a pair such that $h^2=\text{Id}$.
Then there are integers $n_1,n_2$, and $n_3$ such that
$A=n_1A_+\oplus n_2A_-\oplus n_3A_{+-}$.
\end{prop}

\proof Denote by $M_+=\{ a\in M\mid h(a)=a\}$ and by  $M_-=\{ a\in
M\mid h(a)=-a\}$. Consider the subgroup $M'$ of $M$ generated by
the elements $a\in M_+\cup M_-$. Obviously, $M'=M_+\oplus M_-$.
Denote by $p_{+}:M'\to M_{+}$ and $p_{-}:M'\to M_{-}$ the
projections to the factors. Since $M$ is a free $\mathbb
Z$-module, its submodules $M_+$ and $M_-$ are also free $\mathbb
Z$-modules. Let $M_+$ and $M_-$ are freely generated over $\mathbb
Z$, respectively, by $e_1,\dots,e_{k_1}$ and
$e_{k_1+1},\dots,e_{k_1+k_2}$.

Let us show that for any $a\in M$ the element $2a$ belongs to
$M'$. Indeed, for any $a\in M$ we have $a+h(a)\in M_+$, $a-h(a)\in
M_-$ and $2a=(a+h(a))+(a-h(a))$. Therefore all elements of the
abelian group $M/M'$ are of the second order. Note also that, by
definition of the subgroups $M_{\pm}$, if an element $a\in M$ is
such that $2a\in M_+$ (resp. $2a\in M_-$), then $a\in M_+$ (resp.
$a\in M_-$).

We have $M/M'\simeq (\mathbb Z/2\mathbb Z)^{n_3}$ for some integer
$n_3\geq 0$. Let us choose elements $a_1,\dots, a_{n_3}\in M$
whose images generate the group $M/M'$ and let $2a_i=(\alpha
_{i,1},\dots,\alpha_{i,k_1+k_2})$, $i=1,\dots, n_3$, be the
coordinates of $2a_i$ in the base $e_1,\dots, e_{k_1+k_2}$.
Without loss of generality, adding or subtracting elements from
$M'$, we can assume that for each $i=1,\dots, n_3$ each coordinate
$\alpha_{i,j}$ is equal to $0$ or $1$ and there is $j$ for which
$\alpha_{i,j}=1$.

Let us show that the elements $2a_1,\dots ,2a_{n_3}$ are linear
independent over $\mathbb Z$. Indeed, assume that for some
$m_1,\dots,m_{n_3}\in \mathbb Z$ we have $\sum m_i(2a_i)=0$ in
$M'$. Since $M'$ is a free abelian group, we can assume that at
least one of $m_i$, say $m_{n_3}$, is an odd number. In this case
the image of $a_{n_3}$ in $M/M'$ is a linear combination of the
images of the elements $a_1,\dots,a_{n_3-1}$, that is, we obtain a
contradiction.

Let us show that the elements $p_+(2a_1),\dots ,p_+(2a_{n_3})$
(resp. the elements $p_-(2a_1),\dots ,p_-(2a_{n_3})$) are linear
independent over $\mathbb Z$. Indeed, assume that for some
$m_1,\dots,m_{n_3}\in \mathbb Z$ we have $\sum m_ip_+(2a_i)=0$ in
$M_+$. Since $M_+$ is a free abelian group, we can assume that at
least one of $m_i$, say $m_{n_3}$, is an odd number. Then the
element $2a=2\sum m_ia_i\in M_-$ and therefore, $a=\sum m_ia_i\in
M_-$. But it is impossible, since the images of the elements
$a_1,\dots,a_{n_3-1}$ can not generate $M/M'\simeq (\mathbb
Z/2\mathbb Z)^{n_3}$. As a consequence, we obtain that $n_3\leq
\min(k_1,k_2)$.

Let us show that we can choose elements $a_1,\dots, a_{n_3}$,
whose images generate $M/M'$, so that the system of elements
$p_-(2a_1),\dots ,p_-(2a_{n_3})$ can be extended to a free base of
the free $\mathbb Z$-module $M_-$. Indeed, without loss of
generality (if it is needed we renumber the elements
$e_{k_1+1}\dots,e_{k_1+k_2}$), one can assume that
$\alpha_{1,k_1+k_2}=1$. For $i=2,\dots, n_3$ replacing $a_i$ by
$$a'_i=a_i-\alpha_{i,k_1+k_2}(a_1-\sum_{\scriptstyle j=1
}^{k_1+k_2-1}\alpha_{1,j}(1-\alpha_{i,j})e_j),$$ we obtain a new
system $a_1,a'_2,\dots, a'_{n_3}$ such that
\begin{itemize}
\item[($i$)] the images of $a_1,a'_2,\dots, a'_{n_3}$ generate
$M/M'$;
\item[($ii$)] the coordinates of the elements $2a_1,2a'_2,\dots, 2a'_{n_3}$
in the base $e_1,\dots,e_{k_1+k_2}$ are equal to $0$ or $1$;
\item[($iii$)] for each $i=2,\dots,n_3$ the last coordinate of the element $2a'_i$
in the base $e_1,\dots,e_{k_1+k_2}$ is equal to $0$.
\end{itemize}
Now, it is clear that repeating similar replacements $n_3$ times,
we can find new elements $a_1,a_2,\dots, a_{n_3}$ and a new base
$e_{k_1+1},\dots,e_{k_1+k_2}$ of $M_-$ such that
\begin{itemize}
\item[($i$)] the images of $a_1,a_2,\dots, a_{n_3}$ generate
$M/M'$;
\item[($ii$)] for each $i=1,\dots, n_3$ the coordinates $\alpha_{i,j}$
of the element $2a_{i}$ in the new base $e_1,\dots,e_{k_1+k_2}$
are equal to $0$ or $1$;
\item[($iii$)] for each $i=1,\dots,n_3$ the  coordinate $\alpha_{i,j}$
of the element $2a_i$ in the base $e_1,\dots,e_{k_1+k_2}$ is equal
to $0$ for $j>k_1+k_2-i+1$ and $\alpha_{i,k_1+k_2-i+1}=1$.
\end{itemize}

As a result, we can decompose $M_-$ into the direct sum
$M_-=\widetilde M_-\oplus \overline M_-$, where $\overline M_-$ is
generated by $e_{k_1+1},\dots, e_{k_1+k_2-n_3}$ and $\widetilde
M_-$ is generated by $p_-(2a_1),\dots ,p_-(2a_{n_3})$. Denote by
$\widetilde M$ a subgroup of $M$ generated by the elements
$a_1,\dots ,a_{n_3}$ together with the elements
$e_1,\dots,e_{k_1}\in M_+$. It is easy to see that $M=\widetilde
M\oplus \overline M_-$, $\widetilde M$ is invariant under the
action of $h$,  $M_+=\{ a\in \widetilde M\mid h(a)=a\}$, and
$\widetilde M_-=\{ a\in \widetilde M\mid h(a)=-a\}$. Denote by
$\widetilde M'=M_+\oplus \widetilde M_-\subset \widetilde M$ and
by $\widetilde p_+:\widetilde M'\to M_+$ and $\widetilde
p_-:\widetilde M'\to \widetilde M_-$ the projections to the
factors.

Let us show that, by a linear triangular transformation, we can
replace the system of elements $a_1,\dots, a_{n_3}$ chosen above
by $\overline a_1,\dots, \overline a_{n_3}$ so that the images of
$\overline a_1,\dots, \overline a_{n_3}$ generate $\widetilde
M/\widetilde M'$, the elements  $\widetilde p_-(2\overline a_1),$
$\dots ,$ $\widetilde p_-(2\overline a_{n_3})$ generate
$\widetilde M_-$, and the system of elements $\widetilde
p_+(2\overline a_1)$, $\dots $, $\widetilde p_+(2\overline
a_{n_3})$ can be extended to a free base of the free $\mathbb
Z$-module $M_+$. Indeed, as above, without loss of generality (if
it is needed we can renumber the elements $e_{1}\dots,e_{k_1}$),
one can assume that $\alpha_{1,k_1}=1$. Then for $i=2,\dots, n_3$
replacing $a_i$ by
$$a'_i=a_i-\alpha_{i,k_1}(a_1-\sum_{\scriptstyle j=1
}^{k_1-1}\alpha_{1,j}(1-\alpha_{i,j})e_j),$$ we obtain a new
system $a_1,a'_2,\dots, a'_{n_3}$ such that
\begin{itemize}
\item[($i$)] the images of $a_1,a'_2,\dots, a'_{n_3}$ generate
$\widetilde M/\widetilde M'$;
\item[($ii$)] the elements  $\widetilde p_-(2a_1),\widetilde p_-(2a'_2),
\dots ,\widetilde p_-(2a'_{n_3})$ generate $\widetilde M_-$;
\item[($iii$)] the coordinates of the elements $2a_1,2a'_2,\dots, 2a'_{n_3}$ in the
base $e_1,\dots,e_{k_1}, \widetilde p_-(2a_1),\widetilde
p_-(2a'_2), \dots ,\widetilde p_-(2a'_{n_3})$ are equal to $0$ or
$1$;
\item[($iv$)] for $i=2,\dots,n_3$ the coordinate $\alpha _{i,k_1}$
of the element $2a'_i$ in the
base $e_1,\dots,e_{k_1},  \widetilde p_-(2a_1),\widetilde
p_-(2a'_2), \dots ,\widetilde p_-(2a'_{n_3})$ is equal to $0$.
\end{itemize}
It is evident that repeating  similar replacements $n_3$ times, we
can find new elements $\overline a_1,\overline a_2,\dots,
\overline a_{n_3}$ and a new free base $e_{1},\dots,e_{k_1}$ of
the free $\mathbb Z$-module $M_+$ such that
\begin{itemize}
\item[($i$)] the images of $\overline a_1,\overline a_2,\dots, \overline a_{n_3}$ generate
$\widetilde M/\widetilde M'$;
\item[($ii$)] the elements  $\widetilde p_-(2\overline a_1),\widetilde p_-(2\overline a_2),
\dots ,\widetilde p_-(2\overline a_{n_3})$ generate $\widetilde
M_-$
\item[($iii$)] for each $i=1,\dots, n_3$ the coordinates $\alpha_{i,j}$
of the element $2\overline a_{i}$ in the new base
$e_1,\dots,e_{k_1}, \widetilde p_-(2\overline a_1),\widetilde
p_-(2\overline a_2), \dots ,\widetilde p_-(2\overline a_{n_3})$
are equal to $0$ or $1$;
\item[($iv$)] for each $i=1,\dots,n_3$ the  coordinate $\alpha_{i,j}$, $j\leq k_1$,
of the element $2\overline a_i$ in the base $e_1,\dots,e_{k_1},
\widetilde p_-(2\overline a_1),\widetilde p_-(2\overline a_2),
\dots ,\widetilde p_-(2\overline a_{n_3})$ is equal to $0$ for
$j>k_1-i+1$ and $\alpha_{i,k_1-i+1}=1$.
\end{itemize}

As a result, we can decompose $M_+$ into the direct sum
$M_+=\widetilde M_+\oplus \overline M_+$, where $\overline M_+$ is
generated by $e_{1},\dots, e_{k_1-n_3}$ and $\widetilde M_+$ is
generated by $\widetilde p_+(2\overline a_1),\dots ,\widetilde
p_+(2\overline a_{n_3})$. Denote by $\overline M_{+-}$ a subgroup
of $M$ generated by the elements $\overline a_1, h(\overline a_1)
,\dots ,\overline a_{n_3},h(\overline a_{n_3})$. It is easy to see
that $\overline M_{+-}$ is invariant under the action of $h$ and
$\widetilde M_+\oplus \widetilde M_-\subset \overline M_{+-}$.
Therefore $\text{rk}\, \overline M_{+-}=2n_3$ and, by
construction, $(\overline M_{+-},h)\simeq n_3A_{+-}$. Moreover,
$(M,h)=(\overline M_+\oplus \overline M_-\oplus \overline
M_{+-},h)\simeq n_1A_+\oplus n_2A_-\oplus n_3A_{+-}$, where
$n_1=k_1-n_3$ and $n_2=k_2-n_3$. \qed

\begin{lem} \label{char} Let $(M,h)=n_1A_+\oplus n_2A_-\oplus n_3A_{+-}$.
Consider the semi-direct product $G=M\leftthreetimes\langle
h\rangle$. Then \begin{itemize} \item[($i$)] $G/G'\simeq \mathbb
Z^{n_1+n_3+1}\oplus (\mathbb Z/2\mathbb Z)^{n_2}$, \item[($ii$)]
$\det (t\text{Id}-h_{\mathbb C})=(t-1)^{n_1+n_3}(t+1)^{n_2+n_3}$
.\end{itemize}
\end{lem}

\proof The group $G$ can be given by  presentation
\[
\begin{array}{ll} G\simeq  \langle e_1,\dots,e_{n_1+n_2+2n_3},h & \mid
[e_i,e_j]=1,
1\leq i,j\leq n_1+n_2+2n_3, \\
& h^{-\varepsilon}e_ih^{\varepsilon}=e_{n_3+i},\, \, 1\leq i\leq n_3, \varepsilon=\pm 1, \\
& h^{-1}e_{2n_3+i}h=e_{2n_3+i},\, 1\leq i\leq n_1, \\
& h^{-1}e_{2n_3+n_1+i}h=e^{-1}_{2n_3+n_1+i},\, 1\leq i\leq n_2
\rangle .
\end{array}
\]
Therefore the group $G/G'$ can be given by  presentation
\[
\begin{array}{ll} G/G'\simeq  \langle e_1,\dots,e_{n_1+n_2+2n_3+1} & \mid
[e_i,e_j]=1,
1\leq i,j\leq n_1+n_2+2n_3+1, \\
& e_i=e_{n_3+i},\, \, 1\leq i\leq n_3,  \\ &
e_{2n_3+n_1+i}=e^{-1}_{2n_3+n_1+i},\, 1\leq i\leq n_2 \rangle
\end{array}
\]
that is, $G/G'\simeq \mathbb Z^{n_1+n_3+1}\oplus (\mathbb
Z/2\mathbb Z)^{n_2}$.

Part ($ii$) of Lemma \ref{char} is trivial. \qed

Let us return to the proof of Theorem \ref{two}.  Let a polynomial
$\Delta (t)=(-1)^{n+k}(t-1)^n(t+1)^k$ be the Alexander polynomial
of a Hurwitz $C$-group $G$. Then by Theorem 5.9 of \cite{Gr-Ku},
$G$ consists of $n+1$ irreducible components, that is, $G/G'\simeq
\mathbb Z^{n+1}$.

Consider the canonical $C$-epimorphism $\nu: G\to \mathbb F_1$,
$N=\ker \nu$. By \cite{Ku.O}, the group $N$ is finitely presented.
Therefore $N/N'$ is a finitely generated abelian group. Let
$T=\text{Tors}(N/N')$ be the subgroup of $N/N'$ consisting of the
elements of finite orders. Then $T$ is invariant under the action
of $h$. Therefore, $h$ induces an action (we denote it by the same
letter) $h$ on the free abelian group $M=(N/N')/T$. Since the
action of $h_{\mathbb C}$ on $(N/N')\otimes \mathbb C\simeq
M\otimes \mathbb C$ is semi-simple and all roots of the
characteristic polynomial $\det (h_{\mathbb C}-t\text{Id})=\Delta
(t)$ are equal to $\pm 1$, we have $h^2_{\mathbb C}=\text{Id}$.
Consequently, $h^2=\text{Id}$ on $M$. It follows from Proposition
\ref{decomposition} that $(M,h)\simeq n_1A_+\oplus n_2A_-\oplus
n_3A_{+-}$ for some non-negative integers $n_1,n_2,n_3$.

It is easy to see that $T$ is a normal subgroup of $G/N'$.
Therefore the exact sequence
$$1\to
N/N'\to G/N' \buildrel{\nu}\over\longrightarrow \mathbb F_1\to 1
$$ implies the following exact sequence
\begin{equation} \label{exact} 1\to M\to
\overline G\buildrel{\nu}\over\longrightarrow \mathbb F_1\to 1,
\end{equation}
where $\overline G=(G/N')/T$. Denote by $f:G\to \overline G$ and
$g:N\to M$ the canonical epimorphisms, and by $\nu_1:\overline
G\to \mathbb F_1$ the epimorphism induced by $\nu$. Exact sequence
(\ref{exact}) can be included to the following commutative diagram
of exact sequences.
\\

\begin{picture}(0,0)(0,0)
\put(85,-5){$1\longrightarrow  N  \longrightarrow G
\buildrel{\nu}\over\longrightarrow  \mathbb F_1 \longrightarrow
1$} \put(85,-50){$1\longrightarrow  M \longrightarrow \overline G
\buildrel{\nu_{1}}\over\longrightarrow \mathbb F_1 \longrightarrow
1$} \put(121,-10){\vector(0,-1){25}} \put(114,-23){$g$}
\put(156,-10){\vector(0,-1){25}} \put(147,-23){$f$}
\put(192,-10){\vector(0,-1){25}} \put(195,-23){$\simeq$}

\end{picture}
\\
\\
\\
\\
\\

It is easy to see that $\overline G\simeq M\leftthreetimes\langle
h\rangle$. Therefore, by Lemma \ref{char}, we have
\begin{equation} \label{nk}
n=n_1+n_3 \qquad \text{and}\qquad k=n_2+n_3.
\end{equation}
The epimorphism $f$ induces an epimorphism $f_*:G/G'\simeq
Z^{n+1}\to \overline G/\overline G'$. Applying Lemma \ref{char}
one more, we have $\overline G/\overline G'\simeq \mathbb
Z^{n_1+n_3+1}\oplus (\mathbb Z/2\mathbb Z)^{n_2}$. Therefore,
\begin{equation} \label{ineq}
n_1+n_2+n_3+1\leq n+1.
\end{equation}
 It follows from (\ref{nk}) and
(\ref{ineq}) that $n_2=0$ and $k\leq n$.

In \cite{Gr-Ku}, it was shown that the Hurwitz $C$-group
\[
\begin{array}{rl} G(2)=<x_1,\dots,x_4\mid  &
x_2^2x_1x^{-2}_2=x_4,\, x_3=x_2, \, \, x_4^2x_2x_4^{-2}=x_2 \\
 & [x_i,x_1\dots
x_4]=1\,\,\text{for}\, \, i=1,\dots, 4>\end{array}
\]
has the Alexander polynomial $\Delta (t)=(t^2-1)$. Next, the
Alexander polynomial of abelian $C$-group $\mathbb Z^n$ is equal
to $(-1)^{n-1}(t-1)^{n-1}$. Therefore, for $k\leq n$, to prove the
existence of a Hurwitz $C$-group $G$ whose Alexander polynomial
$\Delta (t)=(-1)^{n+k}(t-1)^n(t+1)^k$, it suffices to consider a
Hurwitz product $G=G(2)^{\diamond k}\diamond Z^{n+1-k}$. \qed

\section{Two examples of $C$-groups}
In this section, it will be shown that in  the general case of
$C$-groups, statements similar to Theorem 5.9 of \cite{Gr-Ku} and
Theorem \ref{two}  are not true.
\begin{ex} \label{ex1} There is a $C$-group consisting of two irreducible components
whose Alexander polynomial $\Delta (t)=(t-1)^2$.
\end{ex}

Consider a $C$-group $G$ given by $C$-presentation
$$G=\langle x_1,x_2,x_3\mid x_3=x_1^{-1}x_2x_1,\, \,
x_3=x_1^{-1}x_3x_2x_3^{-1}x_1\rangle .$$

It is easy to see that $G$ is a $C$-group consisting of two
irreducible components. Denote by $N$ the kernel of the canonical
$C$-epimorphism $\nu : G\to \mathbb F_1$. Without loss of
generality, we can assume that $\widetilde h(g)=x_1gx_1^{-1}$ for
$g\in N$.

To find a finite presentation of $N$, let us use the Reidemeister
-- Schreier method. Let us recall briefly it (see particulars, for
example, in \cite{M-K-S}, section 2.3). In the beginning, so
called a Schreier system of representatives should be chosen, that
is, a representative $s_i$ of each right coset of the subgroup $N$
in the group $G=<x\in X\mid r=1, r\in \mathcal R>$ is chosen so
that if a word $s_i$ is a representative, then all its unitial
subwords are also representatives of some cosets (in our case we
choose the elements $x_1^k$, $k\in \mathbb Z$, as Schreier
representatives of the cosets of $N$ in $G$). Then the group $N$
is generated by $a_{i,j}=s_i\cdot x_j\cdot\overline{s_ix_j}^{-1}$,
where $\overline{s_ix_j}$ is the Schreier representative of the
coset containing the element $s_ix_j$, and the system of defining
relations for $N$ is
$$\{ s_irs_i^{-1}=1\mid r\in \mathcal R\} ,$$
where the relations $s_irs_i^{-1}$ are written as words in the
letters $a_{i,j}$ and  $a^{-1}_{i,j}$.

In our case $N$ is generated by the elements
\begin{equation} \label{gen}
a_{k,j}=x_1^kx_jx_1^{-(k+1)}, \end{equation} where $j=2,3$ and
$k\in \mathbb Z$. It is easy to see that the action $\widetilde h$
is given by $\widetilde h(a_{k,j})=a_{k+1,j}$.

The relation $x_3=x_1^{-1}x_2x_1$ gives rise to relations
 \begin{equation} \label{eq2}
a_{k,3} =a_{k-1,2} \end{equation} for $k\in \mathbb Z$, and the
relation $x_3=x_1^{-1}x_3x_2x_3^{-1}x_1$ gives rise to relations
 \begin{equation} \label{eq3}
a_{k,3} =a_{k-1,3}a_{k,2}a_{k,3}^{-1} \end{equation} for $k\in
\mathbb Z$. By the Reidemeister -- Schreier method, the set
consisting of relations (\ref{eq2}) and (\ref{eq3}) is a set of
defining relations for $N$.

Denote by $\overline a_{k,j}$ the image of $a_{k,j}$ in $N/N'$. It
follows from (\ref{eq2}) and (\ref{eq3}) that the abelian group
$N/N'$ is generated by the elements $\overline a_{k,j}$, $j=2,3$
and $k\in \mathbb Z$, being subject to the relations
\begin{equation} \label{eq4}
\begin{array}{lll}
\overline a_{k,3} & = & \overline a_{k-1,2} \\
 2\overline a_{k,3} & = & \overline a_{k-1,3}+ \overline a_{k,2}
 \end{array}
\end{equation} for $k\in \mathbb Z$.
The action of $h$ on $N/N'$, induced by $\widetilde h$, is given
by $h(\overline a_{k,j})=\overline a_{k+1,j}$ for  $j=2,3$ and
$k\in \mathbb Z$.

It follows from (\ref{eq4}) that $N/N'$ is generated by $\overline
a_{0,3}$ and $\overline a_{1,3}$ and we have $\overline
a_{2,3}=-\overline a_{0,3}+2\overline a_{1,3}$. Therefore, in the
base $\overline a_{0,3}$, $\overline a_{1,3}$, the automorphism
$h$ is given by the matrix
$$ h=\left(
\begin{array}{rl}
0 & -1 \\
1 & \, \, \, \, \, 2
\end{array}
\right) .
$$
One can easily check that $\det (h-t\text{Id})=(t-1)^2$.

\begin{ex} \label{ex2} There is a $C$-group consisting of two irreducible components
whose Alexander polynomial $\Delta (t)=(1-t)(t+1)^2$.
\end{ex}

Consider a $C$-group $G$ given by $C$-presentation
$$G=\langle x_1,x_2,x_3\mid x_3=x_1^{-1}x_2x_1,\, \,
[x_1,(x_3^{2}x_1^{-1}x_2x_1^{-1}x_3)]=1 \rangle .$$

Obviously, $G$ is a $C$-group consisting of two irreducible
components. Denote by $N$ the kernel of the canonical
$C$-epimorphism $\nu : G\to \mathbb F_1$. Without loss of
generality, we can assume that $\widetilde h(g)=x_1gx_1^{-1}$ for
$g\in N$. Applying the Reidemeister -- Schreier method, one can
show that $N$ is generated by the elements
\begin{equation} \label{gen}
a_{k,j}=x_1^kx_jx_1^{-(k+1)}, \end{equation} where $j=2,3$ and
$k\in \mathbb Z$.  The action $\widetilde h$ on $N$ is given by
$\widetilde h(a_{k,j})=a_{k+1,j}$.

As in Example \ref{ex1}, the relation $x_3=x_1^{-1}x_2x_1$ gives
rise to the relations
 \begin{equation} \label{eq22}
a_{k,3} =a_{k-1,2} \end{equation} for $k\in \mathbb Z$, and the
relation $[x_1,(x_3^{2}x_1^{-1}x_2x_1^{-1}x_3)]=1$ gives rise to
the relations
 \begin{equation} \label{eq33}
a_{k,3}a_{k+1,3}a_{k+1,2}a_{k+1,3}a_{k+2,3}^{-1}a_{k+2,2}^{-1}a_{k+2,3}^{-1}
a_{k+1,3}^{-1}=1
\end{equation} for $k\in \mathbb Z$. By the Reidemeister --
Schreier method, the set of relations (\ref{eq22}) and
(\ref{eq33}) is a set of defining relations for $N$.

Denote by $\overline a_{k,j}$ the image of $a_{k,j}$ in $N/N'$. It
follows from (\ref{eq22}) and (\ref{eq33}) that the abelian group
$N/N'$ is generated by the elements $\overline a_{k,j}$, $j=2,3$
and $k\in \mathbb Z$, being subject to the relations
\begin{equation} \label{eq44}
\begin{array}{c}
\overline a_{k,3}  =  \overline a_{k-1,2} \\
\overline a_{k,3}+\overline a_{k+1,3}+\overline
a_{k+1,2}-2\overline a_{k+2,3}-\overline a_{k+2,3}=0
 \end{array}
\end{equation} for $k\in \mathbb Z$.
The action of $h$ on $N/N'$, induced by $\widetilde h$, is given
by $h(\overline a_{k,j})=\overline a_{k+1,j}$ for  $j=2,3$ and
$k\in \mathbb Z$.

It follows from (\ref{eq44}) that $N/N'$ is generated by
$\overline a_{0,3}$, $\overline a_{1,3}$ and $\overline a_{2,3}$
and we have $\overline a_{3,3}=\overline a_{0,3}+\overline
a_{1,3}-\overline a_{2,3}$. Therefore, in the base $\overline
a_{0,3}$, $\overline a_{1,3}$, $\overline a_{2,3}$, the
automorphism $h$ is given by the matrix
$$ h=\left(
\begin{array}{rcl}
0 & 0 & \, \, \, \, 1 \\
1 & 0 & \, \, \, \, 1 \\
0 & 1 & -1
\end{array} \right) .
$$
 One can check easily that $$\det (h-t\text{Id})=(1-t)(t+1)^2.$$

 \ifx\undefined\bysame
\newcommand{\bysame}{\leavevmode\hbox to3em{\hrulefill}\,}
\fi

\end{document}